\documentclass[12pt]{article}
\usepackage{amsmath}
\usepackage{amssymb}
\usepackage{amscd}
\usepackage{amsfonts}
\newcommand{\om}{\Omega_N}
\def\reel{\hbox{{\rm R}\kern-1em\hbox{{\rm I} }}}
\def\relatif{\ \hbox{{\rm Z}\kern-.4em\hbox{\rm Z}}}
\def\nat{\hbox{{\rm N}\kern-1em\hbox{{\rm I} } }}
\def\comp{\hbox{{\rm C}\kern-.55em\hbox{{\rm I} } }}
\def\smallcomp{\hbox{\fiverm C}\kern-.35em{\hbox{\fiverm I}}}

\def\fudge{\mathchoice{}{}{\mkern.5mu}{\mkern.8mu}}
\def\bbc#1#2{{\rm \mkern#2mu\vbar\mkern-#2mu#1}}
\def\bbb#1{{\rm I\mkern-3.5mu #1}} \def\bba#1#2{{\rm #1\mkern-#2mu\fudge
#1}}
\def\bb#1{{\count4=`#1 \advance\count4by-64 \ifcase\count4\or\bba
A{11.5}\or \bbb B\or\bbc C{5}\or\bbb D\or\bbb E\or\bbb F \or\bbc
G{5}\or\bbb H\or \bbb I\or\bbc J{3}\or\bbb K\or\bbb L \or\bbb
M\or\bbb N\or\bbc O{5} \or \bbb P\or\bbc Q{5}\orrrr \bb R\or\bbc
S{4.2}\or\bba T{10.5}\or\bbc U{5}\or    \bba V{12}\or\bba
W{16.5}\or\bba X{11}\or\bba Y{11.7}\or\bba Z{7.5}\fi}}
\def\rat{\hbox{{\rm Q}\kern-.70em\hbox{{\rm I} } }}

\def \E {\bbb E}
\newcommand{\vx}{\vert X_N(t)\vert}

\newcommand{\hij}{\phi(i,j)}
\newcommand{\refm [1]} {(\ref{#1})}

\newcounter{corcountrer}
\newcounter{theoremcounter}
\newcounter{lemmacounter}
\newcounter{problemcounter}
\newcounter{remarkcounter}
\newcounter{propositioncounter}

\newtheorem{cor}[corcountrer]{Corollary}
\newtheorem{thm}[theoremcounter]{Theorem}
\newtheorem{lemma}[lemmacounter]{Lemma}
\newtheorem{prob}[problemcounter]{Problem}
\newtheorem{remark}[remarkcounter]{Remark}
\newtheorem{proposition}[propositioncounter]{Proposition}

\newcommand{\wx}{\vert X_N(t)\vert, \ t\ge 0}

\newcommand{\ber}{\begin{eqnarray}}
\newcommand{\ena}{\end{eqnarray}}
\newcommand{\rules}{\ \mbox{\rule{.5mm}{4mm}}\ }

\newcommand{\be}{\begin{equation}}

\newcommand{\nin}{\noindent}
\newcommand{\non}{\nonumber}

\newcommand{\ds}{\displaystyle}
\def\qed{\hfill \vrule height1.3ex width1.2ex depth-0.1ex}
\def\bbb#1{{\rm I\mkern-3.5mu #1}} \def\bba#1#2{{\rm #1\mkern-#2mu\fudge
#1}}

\newcommand{\omr}{\Omega_{N,r}}

\newcommand{\la}{\label}
\newcommand{\en}{\end{equation}}

\newcommand{\PP}{\mathbb{P}}

\title{On time dynamics of coagulation-fragmentation processes
 }

\author{{\bf Boris L. Granovsky}
\thanks{E-mail: mar18aa@techunix.technion.ac.il} \\
Department of Mathematics, Technion-Israel Institute of Technology,\\
Haifa, 32000, Israel \\ \quad {\bf Michael M. Erlihson}
\thanks{E-mail: maerlich@tx.technion.ac.il}\\
Department of Mathematics, Technion-Israel Institute of Technology,\\
Haifa, 32000, Israel.}
\begin{document}
\maketitle \vskip 5cm \nin American Mathematical Society 2000
subject classifications.

\nin Primary-60J27; Secondary-60J35, 05A18, 82C23.

\nin Keywords and phrases:  Time dynamics,
Coagulation-fragmentation models, Gibbs distributions on the set
of integer partitions.
 \newpage
 \pagestyle{myheadings}
 \markboth{ }{\qquad \qquad\qquad\qquad\qquad\qquad\qquad\qquad\qquad transient }

\begin{center}
{\bf Abstract}
\end{center}
We establish a characterization of coagulation-fragmentation
processes, such that the induced birth and death  processes
depicting the total number of groups at time $t\ge 0$ are  time
homogeneous. Based on this, we provide a characterization of
mean-field Gibbs coagulation-fragmentation models, which extends
the one derived by Hendriks et al. As a by- product of our
results, the class of solvable models is widened and a question
posed by N. Berestycki and Pitman is answered, under restriction
to mean-field models.

 \newpage
 \section{Introduction,  objective and the context}

The time dynamics of a time homogeneous Markov process $X(t),\ t\ge 0$
on a space $\Omega=\{\eta\}$ of states $\eta$ is described by the
set of transition probabilities
$$p_{\tilde\zeta}(\eta;t):=\PP(X(t)=\eta\ \mbox{\rule{.5mm}{4mm}}\ X(0)=\tilde\zeta), \ \tilde\zeta,\eta\in
\Omega,\ \ t\ge 0.$$
 Given the rates of
the infinitesimal state transitions, the explicit expressions for the
transition probabilities $p_{\tilde\zeta}$ as solutions of a
Kolmogorov system, are  known only for  a few special cases
of the rates. The corresponding models are called solvable. For
the above reason, time dynamics of Markov processes remain,
generally speaking, a mystery. As an example, even for birth-death
processes on the set of integers, the explicit solutions have been
derived only for a few combinations  of birth and death rates.
This explains why the direction of research in this area turned to the
estimation of the rate of convergence (=spectral gap) of the
transition probabilities as $t\to \infty.$ Nevertheless, hunting
for solvable models continues to be of interest.

In the present paper we pursue the above objective for stochastic
processes of coagulation and fragmentation $(CFP's)$.
  We adopt the formulation of a  $CFP=CFP(N)$ given in
 \cite{DGG} on the basis of classic works of Whittle \cite{W} and Kelly
 \cite{Kel} devoted to deterministic and stochastic models of clustering
 in polymerization, electrical networks and in a variety of other fields.
   A $CFP$ $X_N(t), \ t\ge0$ is defined as a time homogeneous  Markov chain
   on the state space $\om$ of all partitions
   $\eta=(n_1,\ldots,n_N):\sum_{i=1}^Nin_i=N,$ of a given  integer $N$. Here $N$ codes the total population of indistinguishable particles
   partitioned into groups (=clusters) of different  sizes, while $n_i$ is the number
of groups of size $i.$   Infinitesimal (in time) events are a coagulation
   of two groups into one and a fragmentation of one group into two groups, and the basic assumption is that  the rates
   (intensities) of the above two single transitions depend only on sizes of
   groups (and do not depend on $N$).
   Namely, the rate of a single coagulation
   of two groups of sizes $i$ and $j$, such that  $2\le i+j\leq N$, into one group of size $i+j$ is $\psi(i,j)$,
   whereas  the rate of a single fragmentation  of a group with size $i+j$ into two groups of sizes $i$ and $j$ is
   $\phi(i,j)$. The functions $\psi$ and $\phi$ are assumed to be
   non negative and symmetric
   in $i,j$.

  Next, we define the induced rates of infinitesimal state transitions.
   Given a state $\eta\in \Omega_N$ with
  $n_i,n_j>0$ for some $1\le i,j\le N,$ denote by $\eta^{(i,j)}\in\Omega_N$ the state that is obtained
  from $\eta$ by a   coagulation of any two groups  of sizes $i$ and $j$,
  and denote by $K(\eta\rightarrow \eta^{(i,j)})$ the rate of the infinitesimal
  state transition $\eta\rightarrow\eta^{(i,j)}.$ Similarly,
  for a given state  $\eta\in\Omega_N$ with $n_{i+j}>0,$ let
  $\eta_{(i,j)}$ be the state that is obtained from $\eta$ by a fragmentation
  of a group of size $i+j\ge 2$ into two
  groups of sizes $i$ and $j,$ and let $F(\eta\rightarrow\eta_{(i,j)})$ be the rate of the infinitesimal state transition
  $\eta\rightarrow \eta_{(i,j)}.$ We assume that the rate    $K(\eta\rightarrow \eta^{(i,j)})$ is  equal to the sum of rates of all
  single coagulations of  $n_i$ groups  with size $i$
   with $n_j$ groups with size $j $, and   that    $F(\eta\rightarrow\eta_{(i,j)})$
   is the sum
   of rates of all single fragmentations of  $n_{i+j}$ groups with size $i+j$  into
  two groups of sizes $i$ and $j$. As a result, we get the following expressions
  for the rates of state transitions:
  \ber
  \non K(\eta\rightarrow \eta^{(i,j)})&=& n_in_j\psi(i,j),\quad i\neq j,\quad 2\leq
  i+j\leq N ,\\ \non K(\eta\rightarrow \eta^{(i,i)})&=&\frac{ n_i(n_i-1)}{2}\ \psi(i,i),\quad 2\leq
  2i\leq N,\\ F(\eta\rightarrow\eta_{(i,j)})&=& n_{i+j}\phi(i,j),\quad
  2\leq i+j\leq N. \la{3_8}
  \ena
We note that an interpretation of the coagulation kernel $K$ in terms
of the kinetics of droplets of different masses can be found in
\cite{mar}.

  Following \cite{Gran}, we  call $CFP's$ with rates of  state transitions of
  the form \refm[3_8] mean field models,
  meaning  that at any state $\eta\in \Omega_N,$ any group can coagulate with any other one or can be fragmented into
  any two parts. We also note that  (\cite{DGG}) a characterization of positive rates of single transitions
  $\psi(i,j), \hij$ that provide reversibility of mean- field $CFP's$ is known  (\cite{DGG}).

We now describe briefly the context of the present paper. The
paper is devoted to the time evolution of the above mean field
$CFP's$ and it consists of two sections. Section 2 is divided into
three subsections. In Subsection 2.1 we characterize the $CFP's$
$X_N(t),\ t\ge 0$ having  time homogeneous  processes $\vert
X_N(t)\vert,\ t\ge 0$ depicting the total number of groups at time
$t\ge 0$. The key result of the paper, stated precisely in Theorem
1 in Subsection 2.2, establishes the equivalence of the following
two conditions:

(i) The birth and death process $\vert X_N(t)\vert,\ t\ge 0$
 is
time homogeneous;

 (ii) The conditional distribution of a $CFP(N),$
given a total number of groups at time $t\ge 0,$ is a time
independent  Gibbs distribution on the set of partitions of $N$
with given number of components.

Consequently, a  characterization of Gibbs $CFP's,$ which  extends
the one  by Hendriks et al ( \cite{HendSES}), is derived.

In  Subsection 2.3 we discuss the following three topics
related to our main result: Steady state distributions of $CFP's$,
Gibbs $CFP's$ on  set partitions and Spectral gaps of Gibbs
$CFP's.$ In particular, under restriction to mean-field models, we
obtain a negative  answer to a question posed by N. Berestycki and
Pitman (\cite{BerP}) about the existence of certain Gibbs $CFP's$.

\newpage
\section{Main
result}

 We  say that  states $\eta,\tilde{\eta}\in\om$ are neighbors:
$\tilde{\eta}\sim \eta,$ if one of the states is obtained either
by a single coagulation or a single fragmentation of components of
the other state.
   Then the preceding description of a $CFP(N),$ say $X_N^{(\rho)}(t),\ t\ge 0,$  starting from an initial probability
   distribution $\rho$ on $\om$ allows us to write the
    corresponding Kolmogorov system  as follows
\ber \non
\dot{p}_{\rho}(\eta;t)&=&-p_{\rho}(\eta;t)\Big(\sum_{\tilde{\eta}\sim
\eta}\big(K(\eta\rightarrow\tilde{\eta})+F(\eta\rightarrow\tilde{\eta}
)\Big)+ \\ && \sum_{\tilde{\eta}\sim
\eta}p_{\rho}(\tilde{\eta};t)\big(K(\tilde{\eta}\rightarrow\eta)
+F(\tilde{\eta}\rightarrow\eta )\big),\quad
\tilde\eta,\eta\in\om,\quad t\ge 0. \la{kol1} \ena

Note that the seminal system of   Smoluchowski equations (1918)
for pure coagulation  can be viewed as an approximation to
\refm[kol1]
  obtained by neglecting correlations between group numbers at time $t\ge 0$.
This issue is widely discussed in the literature, (see \cite{DGG},
\cite{nor}, \cite{frgr2},\cite{ald1}).

  \subsection{ Process of the total number of groups.}

  In our study of  time dynamics of a $CFP$
$X_N^{(\rho)}(t)=(n_1(t),\ldots,n_N(t))\in \om, \ t\ge 0,$
 a central role is played
by the induced stochastic process

 \ber
  \vert X_N^{(\rho)}(t)\vert:=\sum_{i=1}^{N}{n_i}(t),
  \quad  t\geq0,
  \la{3_3}
  \ena
 which depicts the total number of groups in the generic $CFP$ at time $t\ge 0.$
We denote throughout the paper
$$\Omega_{N,r}=\{\eta\in\om:\vert\eta\vert=r\},\quad  r=1,\ldots, N$$
the set of all partitions of $N$ with exactly $r$ components.

  It follows from the definition of a $CFP(N) $  that $\vert X_N^{(\rho)}(t)\vert, \ t\ge 0$ is
  a Markov  birth and
  death process on  the state space $\{1,2,\ldots, N\},$ with rates of birth and
  death\\
  $\lambda_{r,N},\ 1\leq r\leq N-1,\ \mu_{r,N}, \ 2\leq r\leq N,$ respectively, defined
in a   usual way, as  in \refm[3_500] below.
  However, in
  contrast to the generic $CFP(N),$
  the process $\vert X_N^{(\rho)}(t)\vert, \ t\ge 0$ is, in general,
   not  homogeneous in time, which presents a big problem
for
  the investigation of the process.\\
The following example demonstrates the phenomenon of dependence of
the rates $\lambda_{r,N},\mu_{r,N}$ on time $t\ge 0 $ and on an
initial distribution $\rho$ of the generic $CFP(N),$ that causes
the time-inhomogeneity  of the induced birth and death process.
  \\
{\bf Example:} Consider  a $CFP(N), \ N>4$ of pure coagulation
with $\psi(1,1)=\psi(1,2)=0$ and all other $\psi(i,j)>0.$  It is
clear that  $\Omega_{N,N-2}=\eta_1\bigcup \eta_2,$ where
$\eta_1=(N-3,0,1,0,\ldots, 0),$\\$  \eta_2=(N-4, 2,0,\ldots, 0).$
Assuming that the process starts from an initial distribution
$\rho $ on $\Omega_{N,N-2}$, s.t. $\rho(\eta_i)=p_i>0, \ i=1,2,
\ p_1+p_2=1$ and denoting $A_i=\sum_{\zeta\in
\Omega_{N,N-3}}K(\eta_i\rightarrow \zeta)>0,\ i=1,2,$ we have
$$ \dot{p}_{\rho}(\eta_i;t)=-A_ip_{\rho}(\eta_i;t), \quad t\ge 0,\quad i=1,2,  $$
since the transitions $ \Omega_{N,N-1}\rightarrow \eta_i, \ i=1,2$
are impossible. Hence, $$p_{\rho}(\eta_i;t)=p_ie^{-A_it}, \quad
t\ge 0,\quad i=1,2$$ and consequently, it follows from \refm[shom]
below that the rate of death
$$\mu_{N-2,N}(t;\rho)=\frac{A_1e^{-A_1t}p_1+ A_2e^{-A_2t}p_2}{e^{-A_1t}p_1+
e^{-A_2t}p_2}, \quad t\ge 0$$ depends on $\rho$ and $t$ iff
$A_1=(N-3)\psi(1,3)\neq A_2=\psi(2,2).$ \qed
\vskip .5cm

We further
 assume that the $CFP's$ considered are
   ergodic.

   We first
  distinguish  $CFP's(N)$ which induce  time homogeneous
  processes\\
  $\vert X_N^{( \rho)}(t)\vert, \ t\ge 0,$ under any initial distribution $\rho$ on
   $\om,$ i.e. processes with birth and death rates not
  depending on $t\ge 0$ and  $\rho.$   Let $\lambda_{r,N}(t;\rho),\ 1\leq r\leq N-1$ and
  $\mu_{r,N}(t;\rho),\ 2\le r\le N$ be, respectively,
  the rates of birth and death at time $t\ge 0$,
  of some $\vert X_N^{(\rho)}(t)\vert:$

  \ber
  & &\non \lambda_{r,N}(t;\rho)=\lim_{\Delta t\rightarrow0^+}{\frac{\mathbf{P}
  \Big(\vert X_N^{( \rho)}(t+\Delta
  t)\vert=r+1\rules\ \vert X_N^{(\rho)}(t)\vert=r\Big)}{\Delta t}},\\
& &\mu_{r,N}(t;\rho)=\lim_{\Delta
t\rightarrow0^+}{\frac{\mathbf{P}
  \Big(\vert X_N^{( \rho)}(t+\Delta
  t)\vert=r-1\rules\ \vert X_N^{(\rho)}(t)\vert=r\Big)}{\Delta t}}.
  \la{3_500}
  \ena
 \refm[3_500] tells us that under the assumption of  ergodicity
of the generic $CFP(N),$
the time independence of  birth and death rates implies their independence on $\rho.$

Clearly, the  birth and death rates in \refm[3_500] are implied
respectively, by the rates $\psi$ of single fragmentations and the
rates $\phi$ of single  coagulations of the generic $CFP(N).$
 It turns  out that the required necessary and sufficient condition
of time homogeneity of the process $\vert X_N^{( \rho)}(t)\vert, \ t\ge 0$
has a simple probabilistic meaning.

\begin{lemma}
\label{homlemma}
 $\vert X_N^{(\rho)}(t)\vert, \ t\ge 0$ is a time homogeneous birth and death
 process under any initial distribution $\rho$ on $\om,$
 if and only if the generic $CFP(N)$ is such that for a given $1\le r\le N$
 the sums of rates
 $\sum_{\tilde \eta\sim\eta} K(\eta\rightarrow \tilde\eta),\ \eta\in \Omega_{N,r}$ and
 $\sum_{\tilde \eta\sim\eta} F(\eta\rightarrow \tilde\eta),\ \eta\in \Omega_{N,r}$
 do not depend on $\eta\in \omr.$
  Under the above
condition,  the first sum and the second sum are equal to the rate
of death $\mu_{r,N}$ and to the rate of birth $\lambda_{r,N}$
respectively, so that for any $\eta\in \omr,$ \ber & &\non
\lambda_{r,N}=\ds\lim_{\Delta t\rightarrow0^+}\frac{1}{\Delta t}
 \ \mathbf{P}\Big(\vert X_N^{( \rho)}(t+\Delta
  t)\vert=r+1\rules \  X_N^{( \rho)}(t)= \eta \Big), \ 1\leq r\leq N-1,\\
& &\mu_{r,N}=\ds\lim_{\Delta t\rightarrow0^+}\frac{1}{\Delta t}
 \ \mathbf{P}\Big(\vert X_N^{(\rho)}(t+\Delta
  t)\vert=r-1\rules \  X_N^{(\rho)}(t)= \eta \Big),\ 2\leq r\leq N,
  \la{lm}
  \ena
under any initial distribution $\rho $ on $\Omega_N$ and all $t\ge0.$

\end{lemma}

{\bf Proof}\ \ \  Recalling that $\vert X_N^{(\rho)}(t)\vert, \
t\ge 0$ is Markov, by the Markovian property of the generic
$CFP(N)$, we firstly assume that $\vert X_N^{( \rho)}(t)\vert, \
t\ge 0$ is  time homogeneous, so that for all $\rho $ on $\om$ and
all $t\ge 0,$ $\lambda_{r,N}(t;\rho)=\lambda_{r,N},\ 1\leq r\leq
N-1$ and
  $\mu_{r,N}(t;\rho)=\mu_{r,N},\ 2\le r\le N.$  We now rewrite \refm[3_500] as

\ber& &\non \lambda_{r,N}=\ds\lim_{\Delta
t\rightarrow0^+}{\frac{1}{\Delta t}
  \frac{\sum_{\eta\in \omr}\mathbf{P}\Big(\vert X_N^{(\rho)}(t+\Delta
  t)\vert=r+1\ \rules
  \  X_N^{(\rho)}(t)=\eta
  \Big)\mathbf{P}\Big( X_N^{(\rho)}(t)=\eta\Big)}{\mathbf{P}\Big(\vert X_N^{(\rho)}(t)
  \vert=r\Big)}},\\ & & \non 1\le r\le N-1,\\
& &\non \mu_{r,N}=\ds\lim_{\Delta t\rightarrow0^+}{\frac{1}{\Delta
t}
  \frac{\sum_{\eta\in \omr}\mathbf{P}\Big(\vert X_N^{(\rho)}(t+\Delta t)\vert=r-1\ \rules \
  X_N^{(\rho)}(t)=\eta
  \Big)\mathbf{P}\Big( X_N^{(\rho)}(t)=
  \eta \Big)}{\mathbf{P}\Big(\vert
  X_N^{(\rho)}(t)\vert=r\Big)}},\\
  & & \non  2\le r\le N.\\
\la{lrmr}
  \ena
By the ergodicity and the time homogeneity properties of the
generic $CFP(N),$ the limits \ber && \non f_b(\eta;r,N):=
\ds\lim_{\Delta t\rightarrow0^+}\frac{1}{\Delta t}
 \ \mathbf{P}\Big(\vert X_N^{(\rho)}(t+\Delta
  t)\vert=r+1\ \rules \  X_N^{(\rho)}(t)=\eta \Big), \quad 1\le r\le N-1,\\
& &f_d(\eta;r,N):=\ds\lim_{\Delta t\rightarrow0^+}\frac{1}{\Delta
t}
 \ \mathbf{P}\Big(\vert X_N^{(\rho)}(t+\Delta
  t)\vert=r-1\ \rules \  X_N^{(\rho)}(t)= \eta \Big), \quad 2\le r\le N,
\la{fbd}
  \ena
do not depend on $t\ge 0$ and $\rho,$ for all $\eta\in \omr,$ so that
\ber &
  &\non \lambda_{r,N}=\sum_{\eta\in
 \omr}\frac{f_b(\eta;r,N)\mathbf{P}\Big( X_N^{(\rho)}(t)=
\eta\Big)}{\mathbf{P}\Big(\vert X_N^{(\rho)}(t)\vert=r\Big)},\quad  1\le r\le N-1,\\
& & \mu_{r,N}=\sum_{\eta\in
\omr}\frac{f_d(\eta;r,N)\mathbf{P}\Big( X_N^{(\rho)}(t)=\eta\Big)}
{\mathbf{P}\Big(\vert X_N^{(\rho)}(t)\vert=r\Big)},\quad 2\le r\le
N. \la{shom}\ena

     Next, setting in \refm[shom], $t=0$ and $\rho(\tilde \zeta)=1,$
     for a $\tilde\zeta\in \Omega_{N,r},$ (so that
     $\mathbf{P}\Big( X_N^{(\rho)}(0)=\tilde\zeta\Big)=1$) it is easy to conclude that
     \refm[shom] together with the time homogeneity assumption imply
     $$\lambda_{r,N}=f_b(\tilde\zeta;r,N)=const, \
     \mu_{r,N}=f_d(\tilde\zeta;r,N)=const,$$
  for all $\tilde\zeta\in \omr,$
  which proves the necessity of the condition \refm[lm].
The sufficiency
   of \refm[lm]  follows immediately from \refm[shom], after we observe
that in view of \refm[fbd], the quantities $f_d(\eta;r,N), \
f_b(\eta;r,N)$ are equal respectively, to the sum  of
 rates of single coagulations
 $\sum_{\tilde \eta\sim\eta} K(\eta\rightarrow \tilde\eta)$ and to
 the sum of rates of single fragmentations
 $\sum_{\tilde \eta\sim\eta} F(\eta\rightarrow \tilde\eta)$  at a state
  $\eta\in \omr.$
\qed \vskip .5cm

In the rest of this subseqtion we will treat the case when the
time homogeneity condition in Lemma 1 holds,  writing  simply
$\vert X_N(t)\vert, \ t\ge 0.$  Now our objective will be to
characterize the rates $\psi(i,j), \phi(i,j)$  that provide the
condition \refm[lm].
   The condition \refm[lm] says that for given $N$ and $ r$ each one of the two limits
   in the RHS
  of \refm[lm] is the same for all $\eta\in \omr$ and all $\rho$ on $\Omega_N.$
  Consequently, the above condition
  conforms to
two separate systems of  linear equations, one for $\psi(i,j)$ and
one for $\phi(i,j),$ and  each one consisting of $\vert\omr\vert$
equations for each $1\le r\le N.$ It is easily seen that for a
fixed $N$ there is a variety of solutions to each of these
systems, which are valid for all  $1\le r\le N.$

  For example,  employing the aforementioned meaning of the limits $f_b$ and $f_d,$
one can   verify that for a given $N>3$ the
  following rates depending on $N$ satisfy \refm[lm]:
  \ber
  \psi(i,j)=
    \begin{cases}
     i+j, & \text{if}\quad 2\leq i+j\leq N-1\\
     l_1(N), & \text{if}\quad {i+j=N}
     \end{cases}
   \la{cor1}
   \ena
   and
   \ber
   \phi(i,j)=
  \begin{cases}
    0 , &\ \text{if} \quad  2\le i+j<N\\
    l_2(N) , & \text{if}\quad  i+j=N,
  \end{cases}
   \la{cor2}
   \ena
  where $l_1$ and $l_2$ are arbitrary nonnegative functions.
  \\

However, due to our basic assumption that the rates $\psi$ and
$\phi$ do not depend on $N,$ time homogeneity of the process
$\vert X_N(t)\vert \ t\ge 0$ implies a very special form of the
above rates of singular transitions.
  \begin{proposition}
  \label{homtheorem}
  $\{\wx\}_{N\ge 1}$ is a sequence of  time homogeneous  birth and
  death processes
  induced by a sequence of   $\{CFP(N)\}_{N\ge 1}$ with rates of single transitions
  $\psi(i,j)$ and $\phi(i,j),$
   if and only if
   the above rates are of the form:
  \be
  \psi(i,j)=a(i+j)+b, \ i,j\ge 1, \ a\ge 0,\quad 2a+b\ge 0
  \la{psi}
  \end{equation}
   and
  \be
  v_k:=\sum_{1\le i\le j:\ i+j=k}\phi(i,j)=\phi(1,1)(k-1), \ k\ge 1,
  \la{phi}
  \end{equation}
  where $v(k)$ is the sum of  rates of all possible single fragmentations of a group of size
  $k\ge 2$ into two groups, whereas $v_2=\phi(1,1)\ge 0$ is arbitrary.
  \end{proposition}

{\bf Proof:} We employ the preceding lemma. Assuming that the
processes $\vert X_N(t)\vert,\ t\ge 0$ are  time homogeneous for
all $N\ge 1,$ we apply the second part of \refm[lm] with $r=2$ to
obtain
$$\psi(i,N-i)=\mu_{2,N},\ i=1,\ldots,N-1,\quad N\ge 1.
$$ Therefore,

  \be
  \psi(i,j)=s(i+j), \ i,j\ge 1,
  \la{lll}
  \end{equation}
   where $s$ is some nonnegative  function on integers which are greater or equal to $2.$

  Next, consider the two states $\eta_1,\eta_2\in\Omega_{N,3}, \ N\ge 5:$
  \ber
   \non\eta_1 &=&(2,0,\ldots,0,\overbrace{1}^{n_{N-2}},0,\ldots,0),
   \\  \eta_2  &=&(1,1,0,\ldots,0,\overbrace{1}^{n_{N-3}},0,\ldots,0).
   \la{3_19}
  \ena
 Applying the equation $f_d(\eta_1;3,N)=f_d(\eta_2;3,N),$ gives
 \ber
 2\psi(1,N-2)+\psi(1,1)=\psi(1,N-3)+ \psi(N-3,2)+\psi(1,2),
 \la{psi_eq}
 \ena
   which by virtue of \refm[lll], is equivalent to $$2s(N-1)+s(2)=s(N-2)+s(N-1)+s(3), \ N\ge 5.$$ Taking into account
  that the last relation should hold for all $N\ge 5,$ we rewrite it as $s(k)-s(k-1)=s(3)-s(2), \ k\ge 3,$
  which proves the necessity of \refm[psi]. For the proof of the necessity of \refm[phi] we consider the quantities
  $f_b(\eta;2,N)$ for $N$ fixed and all states $\eta$ of the form
 $$\eta=(0,\ldots,0,\ldots, 0,\overbrace{1}^{i},0,\ldots, 0, \overbrace{1}^{N-i},0,\ldots,0)\in\Omega_{N,2},
  \quad 1\le i\le N-1.$$  Using   the notation in \refm[phi],
  the condition that $f_b(\eta;2,N)$ should be the same for all the above
  $\eta$ can be written as
  \ber
  v(i)+v(N-i)=const,\  1\le i\le N-1,
  \la{phi_eq}
  \ena
 or, equivalently, $v(N-1)-v(N-2)=v(2)- v(1)=v(2).$
 Since the latter relationship should hold for all $N\ge 2,$ it implies
 \refm[phi]. We turn now to the proof of sufficiency of the conditions
 \refm[psi] and \refm[phi]. Supposing that \refm[psi]
 holds, we have for a state $\eta\in\omr:$
 \ber
 \non f_d(\eta;r,N)&=&\sum_{1\le i<j\le N}\psi(i,j)n_in_j+
  \sum_{1\le i\le N}\psi(i,i)\frac{n_i(n_i-1)}{2}\\
 \non &=& {1\over 2}\Big(\sum_{1\le i,j\le N}\psi(i,j)n_in_j-\sum_{1\le
 i\le N}\psi(i,i)n_i\Big)\\
  \non &=& {1\over 2}\Big(\sum_{1\le i,j\le N}\Big(a(i+j)+ b\Big)n_in_j-\sum_{1\le i\le
  N}\Big(2i a+b\Big)n_i\Big)\\
 \non &= &{1\over 2}\Big(2aNr+br^2-2aN-br\Big), \quad r=2,\ldots, N,\\
 \non f_b(\eta;r,N)&=&\sum_{1\le k\le N} v(k)n_k\\
 &=&\sum_{1\le k\le N} v(2)(k-1)n_k=v(2)(N-r),\quad
 r=1,\ldots,N-1.
 \la{yui}
 \ena
\qed

\begin{cor}
\label{rat}
 The rates of  death and birth
  of a time homogeneous Markov process\\
 $\wx$
 are given by
  \ber
  \non\mu_{r,N}&=&{(r-1)\over 2}\Big(2aN+rb\Big),\quad 2\leq r\leq N,\\
  \lambda_{r,N}&=&\phi(1,1)(N-r),\quad 1\leq r\leq N-1.
  \la{cor4}
  \ena
\end{cor}

\begin{remark}
  \la{SolvRem1}
  \quad\\ (i) The  birth and death process $\wx$
  with rates given by \refm[cor4], has the following interpretation, not related to the  generic $CFP(N).$ Consider a
  nearest neighbor spin system (for reference see \cite{Lig}) of \ $"0"$-s
  and \ $"1"$-s on a complete graph on $N$ vertices
  (sites). Assume that one of the  sites is occupied with a
  $"1"$ which never flips, while spins at all other sites perform flips $0\rightarrow 1$ and $1\rightarrow
  0$ with rates $\tilde{\lambda}_{r,N} $ and $\tilde{\mu}_{r-1,N}$ respectively, where
  $r$ is the total number of  sites of the graph
   occupied by $"1"$-s. (The latter says that a site occupied by a $"1"$ has $r-1$ neighbors occupied by $1$-s
  and a site occupied by a $"0"$ has $r$ such neighbors).
  Consequently, at a state with
  $r\ge 1$ $"1"$-s, the total rate of \ $0\rightarrow 1$ flips is
  $\lambda_{r,N}:=(N-r)\tilde{\lambda}_{r,N}$ and  the total
  rate of\  $1\rightarrow 0$ flips is $\mu_{r,N}:=(r-1)\tilde{\mu}_{r-1,N}.$
  Therefore, the induced birth and death process, say $\zeta_N(t),\ t\ge 0,$ on
$\{1,\ldots,N\}$
  depicting
  the number of sites occupied by $"1"$-s at time $t\ge 0$ is
  Markov and time homogeneous.
   Clearly, if
  $$\tilde{\lambda}_{r,N}=\phi(1,1),\quad 1\le r\le N-1,$$
  $$\tilde{\mu}_{r-1,N}=\frac{1}{2}
  \Big(2aN+br)\Big),\quad 2\le r\le N,$$ the process $\zeta_N(t)$ conforms to  the
  process $\wx,$ associated with  CFP's(N) given by \refm[psi],\refm[phi]. Finally, it is appropriate to note
  that after interchanging the roles of $"0"$-s and $"1"$-s, the spin system with the rates $\tilde{\lambda}_r,$
  $\tilde{\mu}_{r-1}$ as above, is known (for   $N$ fixed) as   a contact
  process.\\
  (ii) It follows from Proposition~\ref{homtheorem} that the
   class of CFP's(N) that induce  time homogeneous processes $\wx $ includes
  processes of pure coagulation ($\phi(1,1)=0$ in \refm[phi]) and
  processes of pure fragmentation ($a=b=0 $ in \refm[psi]).
Also observe that the time homogeneity requirement determines
uniquely the form of  rates of singular coagulations, while it
leaves a certain freedom in the choice of rates of singular
fragmentations.

\end{remark}

  As far as we know,  there are no explicit solutions, i.e. explicit formulae for transition
  probabilities $\PP(\vx=r), \  t\ge 0, \ 1\le r\le N,$ for  birth-death processes
   with the rates given by
  \refm[cor4], when $a,b>0,\ \psi(1,1)>0$ and the initial distribution is
   concentrated on some state $\zeta\in \om.$ The problem here is that the birth and death rates
  in \refm[cor4] are polynomials in $r$ of different degrees, which are 1 and 2 respectively.
  A survey of solvable  birth-death processes with polynomial rates is given in
   \cite{Val}.\\
We will see in the next subsection that under the above condition
of time homogeneity of $\vert X_N(t)\vert, \ t\ge 0$ and certain
initial distributions $\rho,$ the corresponding $CFP's(N)$ are
solvable.
\subsection{Solvable CFP's}

Let a $CFP(N)$
 considered start from an initial distribution
$\rho$ on $\Omega_{N}, $ with  projections $\rho_r$ on the sets
$\omr, \ r=1,\dots,N$: \be \rho_r(\eta):=
  \begin{cases}
    \rho(\eta\ \rules \ \vert \eta
\vert=r),  \quad \eta \in \omr , & \text{if}\quad
\rho(\omr):=\rho(\vert \eta \vert=r)>0\\
    0,\quad \eta \in \omr \ , & \text{if}\quad  \rho(\omr)=0.
  \end{cases}
\end{equation}

It is in order to note that the set of all distributions $\rho$ on
$\om$ with given projections  $\rho_r,$ is \be
\{\rho:\rho(\eta)=\rho_r(\eta)\rho(\omr), \quad \eta\in\omr,\quad
\sum_{r=1}^N\rho(\omr)=1, \quad \rho(\omr)\ge 0,\quad
r=1,\ldots,N\},\la{durk}
\end{equation}
i.e.  the projections $\rho_r, \ r=1,\ldots,N$ define the
associated distribution $\rho $ up to the factors $\rho(\omr),\
r=1,\ldots,N$ in \refm[durk].

Accordingly, we write
 \ber \non p_{\rho}(\eta;t)&=&\PP\Big(X_N^{(\rho)}(t)=\eta\ \rules \ \vert
X_N^{(\rho)}(t)\vert=\vert\eta\vert\Big) \ \PP\Big(\vert
X_N^{(\rho)}(t)\vert=\vert\eta\vert\Big)\\
&:=&  Q(\eta,\rho;t)\ b(\vert\eta\vert, \rho;t), \ \eta\in\om,\
t\ge  0,\quad
 \la{prob1} \ena where  $Q(\eta,\rho;t)$ and $b(\vert\eta\vert,\rho;t)$ denote respectively
the first and the second factors in the RHS of the first  line,
while the conditional probability $Q$ obeys the  initial
conditions \be Q(\eta,\rho;0)= \rho_r(\eta),\ \eta\in\Omega_{N,r}\
, \quad r=1,\ldots,N, \la{in1}
 \end{equation}
 for any $\rho$ on $\om,$ which follow  from the
definitions of $Q$ and $\rho_r.$

We will be interested in $CFP's$ possessing a conditional
probability $Q$ not depending of time under certain initial
distributions $\rho$. If this is the case, it follows from
\refm[in1] that \be Q(\eta,\rho;t)=\rho_r(\eta), \ \eta\in \omr, \
r=1,\ldots,N, \ t\ge 0 \la{toc}.\en Obviously, time independence
\refm[toc] holds for any $CFP(N)$ starting from its stationary
distribution. In view of this, we adopt the following convention.

{\bf Definition 1}
 A $CFP(N)$  possesses a time independent conditional probability
$Q$ if \refm[toc] holds for
 certain projections $\rho_r$ on $\omr, \ r=1,\ldots,N$ and all initial
distributions $\rho$ on $\om$ from the associated set \refm[durk].

 Next we write  the Kolmogorov system
for a  birth and death process $\vert X_N^{(\rho)}(t)\vert, $ $
t\ge 0 $\\ with rates $\lambda_{r,N}(t;\rho),\quad
\mu_{r,N}(t;\rho):$ \ber \non
\dot{b}(r,\rho;t)&&=-b(r,\rho;t)\big(\lambda_{r,N}(t;\rho)+\mu_{r,N}(t;\rho)\big)+
b(r+1,\rho;t)\mu_{r+1,N}(t;\rho)+\\ && b(r-1,\rho;t)
\lambda_{r-1,N}(t;\rho),  \quad r=1,\ldots,N, \la{kol3}\ena where
$\quad b(0,\rho;t)=b(N+1,\rho;t)=0, \quad t\ge 0.$

The following assertion is crucial for our study.

 \begin{proposition}
 \label{ab}

The following two conditions (i) and (ii) are equivalent.

(i) A $CFP(N)$ possesses a conditional probability $Q$ independent
of time $t\ge 0;$

(ii) The  birth and death process $\vert X_N(t),\vert,\ t\ge 0$ is
time homogeneous.

Moreover, the projections $\rho_r,\ r=1,\ldots,N$ defining by \refm[toc] the
time independent conditional probability $Q$  are the unique solution
of the two systems of \
 equations: \ber \mu_{r+1,N}\rho_r(\eta)&=&
\sum_{\zeta\in \ \Omega_{N,r+1}:\ \zeta\sim
\eta}\rho_{r+1}(\zeta)K(\zeta\rightarrow
\eta),\quad \eta\in\Omega_{N,r},\quad r=1,\ldots,N-1, \la{eq1}\\
\lambda_{r,N}\rho_{r+1}(\zeta)&=& \sum_{\eta\in \ \Omega_{N,r}:\
\eta\sim\zeta}\rho_r(\eta)F(\eta\rightarrow \zeta),\quad \zeta\in
\ \Omega_{N,r+1},\quad r=1,\ldots,N-1,
 \la{eq2}
\ena where the rates of state transitions   $F$ and $K$ are given
by \refm[3_8],\refm[psi], \refm[phi], while the rates of birth and
death are as in \refm[cor4].

\end{proposition}

{\bf Proof:} We substitute \refm[prob1] in the Kolmogorov system
\refm[kol1] to obtain
\ber \non &&\dot{b}(r,\rho;t)Q(\eta,\rho;t)+\dot{Q}(\eta,\rho;t)b(r,\rho;t)=\\
\non && -Q(\eta,\rho;t)b(r,\rho;t)\Big(\sum_{\zeta\in \
\Omega_{N,r-1}:\ \zeta\sim \eta}
K(\eta\rightarrow\zeta)+\sum_{\zeta\in \ \Omega_{N,r+1}:\
\zeta\sim \eta}F(\eta\rightarrow\zeta )\Big)+\\\non &&
b(r+1,\rho;t) \sum_{\zeta\in \Omega_{N,r+1}:\ \zeta\sim
\eta}Q(\zeta,\rho;t)K(\zeta\rightarrow\eta)+
b(r-1,\rho;t)\sum_{\zeta\in \Omega_{N,r-1}:\ \zeta\sim
\eta}Q(\zeta,\rho;t)F(\zeta\rightarrow\eta ),\\ \non &&
\eta\in\Omega_{N,r},\quad r=1,\ldots, N,\quad t\ge 0,\\ &&
\Omega_{0,N}=\Omega_{N,N+1}=\varnothing, \quad
b(0,\rho;t)=b(N+1,\rho;t)=0, \quad t\ge 0.
 \la{kolm} \ena

 We firstly prove the implication $(ii)\Rightarrow (i).$  We
substitute in the LHS of \refm[kolm] the expression for
the derivative $\dot{b}(r,\rho;t)$ from \refm[kol3], assuming that
  the process $\vert X_N(t)\vert, \ t\ge 0$ is time
homogeneous and that the initial distribution is $\rho$. Then, by
virtue of Lemma~\ref{homlemma}, the system \refm[kolm] becomes
\ber \non &&
Q(\eta,\rho;t)\Big(b(r+1,\rho;t)\mu_{r+1,N}+b(r-1,\rho;t)
\lambda_{r-1,N}\Big)+\dot{Q}(\eta,\rho;t)b(r,\rho;t)= \\
\non && b(r+1,\rho;t) \sum_{\zeta\in \Omega_{N,r+1}:\ \zeta\sim
\eta}Q(\zeta,\rho;t)K(\zeta\rightarrow\eta)+
b(r-1,\rho;t)\sum_{\zeta\in \Omega_{N,r-1}:\ \zeta\sim
\eta}Q(\zeta,\rho;t)F(\zeta\rightarrow\eta ),\\ \non &&
\eta\in\Omega_{N,r},\quad r=1,\ldots, N,\quad t\ge 0,\\ &&
\Omega_{0,N}=\Omega_{N,N+1}=\varnothing, \quad
b(0,\rho;t)=b(N+1,\rho;t)=0, \quad t\ge 0,\  \text{for\ \ all}\ \
\rho \ \text{on}\  \om,
 \la{kolm3} \ena
where by the   assumption made, the rates
$K(\zeta\rightarrow\eta)$ and $F(\eta\rightarrow\zeta)$ of state
transitions are implied by \refm[psi], \refm[phi] respectively and
the birth and death rates are as in Corollary 1.  Given an initial
distribution $\rho,$ constants $a,b$ and fragmentation rates
$\phi(i,j)$ obeying \refm[phi], a finite Kolmogorov system
\refm[kolm3] has a unique solution $Q,$ provided $a^2+b^2+
\phi(1,1)>0.$ In particular, \refm[kolm3] is satisfied by the
time- independent $Q ,$ such that \refm[toc] holds for all initial
distributions $\rho$ with the projections $\rho_r,\ r=1,\ldots,N$
that obey \refm[eq1],\refm[eq2].

To prove the  implication $(i)\Rightarrow (ii)$, we observe  that
by virtue of \refm[shom], the condition $(i)$ implies that the
birth and death rates do not depend on $t\ge 0.$ By our remark
after \refm[3_500] this leads to the conclusion that the rates do
not depend on $\rho $ either.

Next,   we set $t=0$ in \refm[kolm3] to derive  by virtue of
Definition 1 and the time homogeneity of $\vert X_N(t)\vert,\ t\ge
0,$ that the projections $\rho_r, \ r=1,\ldots,N$ should obey
\refm[eq1] and \refm[eq2].

  It is left to show  the existence and
uniqueness of the solution $\rho_r,\  r=1,\ldots,N$ for the system
of equations \refm[eq1], \refm[eq2], where the rates of state
transitions are induced by $\psi$ and $\phi$ as in
\refm[psi],\refm[phi]. Recalling Lemma~\ref{homlemma}, we treat
the ratios
$$P_C(\zeta\rightarrow \eta):=\frac{K(\zeta\rightarrow
\eta)}{\mu_{r+1,N}},\quad \zeta\in\Omega_{N,r+1},\  \eta\in
\Omega_{N,r},\ \zeta\sim\eta,
$$
$$ r=1,\ldots,N-1$$
as the one- step transition probabilities of a discrete time
nearest-neighbor "coagulation" random walk on the set of
partitions $\om.$ Then $\rho_r(\eta),\ \eta\in \omr$, in the  set
of equations \refm[eq1] can be interpreted as the probability that
the random walk starting at $\eta^*=(N,0,\ldots,0)\in
\Omega_{N,N}$ reaches a given  state $\eta\in \Omega_{N,r} $ at
the $(N-r)-th$  step, so that $\zeta^*=(0,\ldots,1)$ is the
absorbing state. In a similar manner, we consider the nearest
neighbor "fragmentation" random walk on $\om$ with the transition
probabilities
$$P_F(\eta\rightarrow
\zeta):=\frac{F(\eta\rightarrow \zeta)}{\lambda_{r,N}},\quad
\eta\in\Omega_{N,r},\  \zeta\in \Omega_{N,r+1},\ \zeta\sim \eta,
$$
$$ r=1,\ldots,N-1,$$ that starts at $\zeta^*=(0,\ldots,0,1)\in \Omega_{N,1}.$
In this case, $\rho_r(\eta),\  \eta\in\Omega_{N,r} $ in the set of
equations \refm[eq2] is the probability that the "fragmentation"
random walk reaches a given state $\eta\in\Omega_{N,r} $ at the
$(r-1)-th$ step, $\eta^*=(N,\ldots,0)$ being the absorbing state.
Clearly, each one of the two systems has a unique solution
$\rho_r,\ r=1,\ldots, N$ whenever $a^2+b^2>0$ in the first case
and $\phi(1,1)>0$ in the second case.

We demonstrate that when $(a^2+b^2)\phi(1,1)>0$ (=both coagulation
and fragmentation hold), the two systems of equations have the
same  solution  if and only if the transition probabilities $P_C$
and $P_F$ are related in the following way. Let $\rho_r, \
r=1,\ldots,N$ be the probabilities corresponding to the
"coagulation" random walk, under some fixed $a,b:a^2+b^2>0$ in
\refm[psi].
 Then the equations \refm[eq2]
for the "fragmentation" random walk have the same solution
$\rho_r, \ r=1,\ldots,N$ if and only if \be
\rho_r(\eta)P_F(\eta\rightarrow \zeta)=\rho_{r+1}(\zeta)
P_C(\zeta\rightarrow \eta),\ \eta\in\Omega_{N,r}, \ \zeta\in
\Omega_{N,r+1},\ \eta\sim\zeta. \la{F}\end{equation} The
sufficiency of \refm[F] is seen immediately, while the necessity
can be  derived from the following general reasoning, based on the
observation that each one of
 the equations \refm[eq1] and \refm[eq2] is time reversal of the other one.
 Let $\rho_r,\ r=1,\ldots, N$ be the common solution of \refm[eq1] and \refm[eq2].
Then from \refm[eq1], applied for $r=N-2,$ we conclude that under
given $\rho_{N-1}$ and $ \rho_{N-2},$ the values $a,b$ in
\refm[psi] are uniquely determined. Hence,  $\rho_r,\
r=1,\ldots,N$ uniquely determine all  probabilities $P_C$ in
\refm[eq1] induced by \refm[psi]. If now some $\tilde{P}_F$
satisfies \refm[eq2] under the above $\rho_r,\ r=1,\ldots,N,$
  then \refm[eq1] should be satisfied
by $$\tilde{P}_C(\zeta\rightarrow
\eta)=\frac{\rho_r(\eta)\tilde{P}_F(\eta\rightarrow
\zeta)}{\rho_{r+1}(\zeta)}.$$ The aforementioned uniqueness of the
probability $P_C$ proves the claim. (In the discussion following
the proof we find explicitly the solution $\rho_r, \ r=1,\ldots,N$
and demonstrate that the rates of singular fragmentations derived
from \refm[F] satisfy the condition \refm[phi]). \qed

 \vskip.5cm
 Our next purpose will be to find explicitly the
solution $\rho_r(\eta),\ \eta\in\omr, \ r=1,\ldots,N$ of
\refm[eq1],\refm[eq2], in the case when  the time homogeneous
process $\vert X_N(t)\vert,\ t\ge 0$ is given by \refm[cor4].
Since the sets $\Omega_{N,1},\Omega_{N,N}$ are singletons, it
follows from the definition of the conditional probability
$\rho_r$ that $\rho_N(\eta^*)=1,\quad \rho_1(\zeta^*)=1.$
  The
following two cases should be broadly distinguished.

{\bf Case 1:} Non zero coagulation, i.e.
 $a^2+b^2>0.$\\ Following the illuminating
idea of Hendriks et al (\cite{HendSES}), we will seek the
probabilities $\rho_r$ in question in the form

$$  \rho_r(\eta)=\rho_{N,r}(\eta)=\Big(B_{N,r}\Big)^{-1}\
\frac{a_1^{n_1}
a_2^{n_2} \ldots a_N^{n_N}} {n_1!n_2!  \ldots n_N!},$$
\be
\eta=(n_1,\ldots,n_N)\in \Omega_{N,r},\quad r=1,\ldots, N,
\la{Q}
\end{equation} where $B_{N,r}$ is the normalizing constant (= partition
function) known as the $(N,r)$ partial Bell polynomial (see e.g.
\cite{BerP}, \cite{Pi}) induced by the sequence of weights
$\{a_k\}_1^\infty$ that do not depend neither on $N$ nor $r$.
 It follows from \refm[Q] that for given
$\eta=(n_1,\ldots,n_N)\in \Omega_{N,r},$ such that $n_{i+j}>0$ for
some $2\le i+j\le N,$ and $\zeta=\eta_{(i,j)}\in\Omega_{N,r+1},$

\ber
\frac{\rho_{r+1}(\zeta)}{\rho_r(\eta)}=\Big(\frac{B_{N,r}}{B_{N,r+1}}\Big)
  \begin{cases}
    \Big(\frac{a_ia_j}{a_{i+j}}\Big)\
\Big(\frac{n_{i+j}}{(n_i+1)(n_j+1)}\Big), & \text{if}\ i\neq j \\
\Big(\frac{a_i^2}{a_{2i}}\Big)\
\Big(\frac{n_{2i}}{(n_i+1)(n_i+2)}\Big)     , & \text{if}\ i=j.
  \end{cases}
 \la{rat}\ena

Hence, setting, in accordance with Proposition~\ref{homtheorem}
and \refm[3_8], \be K(\eta_{(i,j)}\rightarrow \eta)=
  \begin{cases}
 (a(i+j)+b)   (n_i+1)(n_j+1) , & \text{if}\ i\neq j \\
  (2ia+b)  \frac{(n_i+1)(n_i+2)}{2} , & \text{otherwise},
  \end{cases} \la{FR}
\end{equation} where $a^2+b^2>0$, the  equations \refm[eq1] conform to

\ber \non
\mu_{r+1,N}&=&\Big(\frac{B_{N,r}}{B_{N,r+1}}\Big)\sum_{k=2}^N\frac{(ak+b)\sum_{i+j=k}
a_ia_j}{2a_k}\ n_k,\\ &&(n_1,\ldots,n_N)\in \Omega_{N,r}, \quad
r=1,\ldots,N-1. \la{kash}\ena

Since the RHS of \refm[kash] should not depend on $\eta\in
\Omega_{N,r}$ the equations are solved by the weights defined
recursively by \ber a_k=\frac{(ak+b)\sum_{i+j=k} a_ia_j}{2(k-1)},
\quad k\ge 2, \quad a_1=1.\la{rec}\ena

This is just the solution obtained, by  quite different
considerations,  in \cite{HendSES} (see (18) there), for pure
coagulation processes.

Continuing \refm[kash], we get \ber
\mu_{r+1,N}&=&\Big(\frac{B_{N,r}}{B_{N,r+1}}\Big)\sum_{k=2}^N(k-1) n_k,\\
&&(n_1,\ldots,n_N)\in \Omega_{N,r}, \quad r=1,\ldots,N-1,
\la{kash1}\ena which leads to
 the following relation between
the constants $\mu_{{r+1},N},\  B_{N,r},\ B_{N,r+1} $ induced by
the weights \refm[rec]: \be
\mu_{r+1,N}=(N-r)\Big(\frac{B_{N,r}}{B_{N,r+1}}\Big), \quad
r=1,\ldots, N-1. \la{const}
\end{equation}
Taking into account that $B_{N,N}=\frac{a_1^N}{N!}=(N!)^{-1},$ we
get the explicit expressions for the Bell polynomials in the case
considered: \be B_{N,r}=\frac{\prod_{l=r+1}^N\mu_{l,N}}{N!(N-r)!},
\quad r=1,\ldots,N-1,\la{Bell}
\end{equation}
where $ \mu_{l,N}$ as in \refm[cor4]. Remarkably, the expression
\refm[Bell] for the Bell polynomials enables us to find explicitly
the weights $a_k, \ k\ge 1,$ without solving the recurrence
relation \refm[rec]. In fact, by \refm[Q],

\be a_N=B_{N,1}=\frac{\prod_{r=2}^N\mu_{r,N}}{N!(N-1)!}, \quad
N\ge 2,\end{equation} which can be written as
 \be a_1=1,\quad
a_k=\frac{\prod_{r=2}^k\big(ka+\frac{br}{2}\big)}{k!},\quad k\ge
2, \ 2a+b>0,\ a\ge 0. \la{ak}
\end{equation}

\begin{remark}

The recurrence relation \refm[rec] can be viewed as a modification
of the classic  convolution formula,
$$a_k=\frac{1}{2}\sum_{i+j=k} a_ia_j, \quad k=2,\ldots,\quad
a_1=1,$$ which determines the Catalan numbers (see e.g.
\cite{Lan}). It is interesting to find   the generating function
$g(x)=\sum_{k=1}^\infty a_kx^k$ for the sequence of weights
$\{a_k\}_1^\infty,$ defined by \refm[rec]. Setting
$y(x)=\frac{g(x)}{x},$ it follows from \refm[rec] that the
function $y$ obeys the differential equation

$$y^\prime(1-axy)=y^2a_2,\quad a_2=\frac{2a+b}{2},\quad y(0)=a_1=1,$$
which implicit solution is given by
$$y(x)=\Big(1+\frac{b}{2}xy\Big)^{\frac{2a+b}{b}},\ b>0.$$
\end{remark}

 We now  recover the
fragmentation rates given by \refm[F] in the case of coagulation
rates \refm[FR]. By \refm[rat],\refm[const] and \refm[cor4] we
have \be F(\eta\rightarrow \zeta)=\phi (1,1)\
  \begin{cases}
   \Big(\frac{a_ia_j\
n_{i+j}}{a_{i+j}}\Big)\ \Big(a(i+j)+b\Big) , & \text{if}\ i\neq j \\
    \Big(\frac{a_i^2\
n_{2i}}{2a_{2i}}\Big)\ \Big(2ai+b\Big) , & \text{if}\ i=j.
  \end{cases}
\la{fragm}
\end{equation}
for $\eta=(n_1,\ldots,n_N)\in \Omega_{N,r},$ \ such that
$n_{i+j}>0$\ \ for some $2\le i+j\le N,$ and
$\zeta=\eta_{(i,j)}\in\Omega_{N,r+1}.$ Note that by \refm[rec],
the rates of singular fragmentations induced by  \refm[fragm]
satisfy the condition \refm[phi]. This latter condition appears to
have a physical meaning in the context of $CFP's$ describing
polymerization (see \cite{VAN}).

Also note  that in the case considered, $\rho_r(\eta)>0, \
\eta\in\omr, \ r=1,\ldots,N$ and that the rates of single
transitions are determined by \refm[fragm] up to the constants
$a,b$ and $\phi (1,1).$

{\bf Case 2:} Pure fragmentation.\\ It is clear from the preceding
discussion that under the fragmentation rates of the form
\refm[fragm] (with $a,b$ that are not related to coagulation
rates) the solution $\rho_r$ of \refm[eq2] is given by \refm[Q].
However, in contrast to the case of pure coagulation,
Proposition~\ref{homtheorem} leaves freedom for the choice of
rates of single fragmentations obeying \refm[phi]. In view of
this,  the probabilities $\rho_r$ solving \refm[eq2] will depend
on a particular choice of the above rates, so that $\rho_r$, will
be of Hendriks et al  form \refm[Q] if and only if the  rates of
single fragmentations  are induced by  \refm[fragm]. This is
illustrated by the toy example below. (We recall that under all
above choices of rates of single fragmentations, the rates of the
induced pure birth process remain the same: $\lambda_r=\phi
(1,1)(N-r),\ r=1,\ldots, N$).

 {\bf Example }
Let \ber \phi(i,j)= \phi (1,1)\begin{cases} (i+j-1), & \text{if}\
i=1\ \text{or}\ j=1\\ 0, & \text{otherwise},
\end{cases}
 \ena
 so that \refm[phi] holds.
The corresponding fragmentation random walk is in effect a
deterministic chain on $N$ states $\zeta_1,\ldots, \zeta_r,\ldots,
\zeta_N,$ such that
  $$\zeta_r =(r-1,0,\ldots,0,\overbrace{1}^{N-r+1},0,\ldots,0)\in\omr, \quad 1\le r\le N-1,$$
$$\zeta_N=(N,0,\ldots,0).$$
Consequently, \refm[eq2] implies $\rho_r(\eta)={\bf
{1}}_{\zeta_r}(\eta), \ \eta\in\omr,$ which is not of the form
\refm[Q].

The preceding discussion is summarized in our main theorem.

  \begin{thm}
  \label{mainheorem}{\bf Solvable $CFP's$}\\
Mean -field CFP's $X_N^{(\rho)}(t), \ t\ge 0$ with rates of single
transitions \refm[psi] and \refm[phi] and initial distributions
$\rho$ on $\om,$ s.t. $\rho_r$ on $\omr, \ r=1,\ldots,N$ satisfy
\refm[eq1],\refm[eq2], have time dynamics given by \ber
p_{\rho}(\eta;t)=\rho_r(\eta)\ b(r,\rho;t), \
\eta\in\Omega_{N,r},\ \ r=1,\ldots, N,\ t\ge 0, \la{dyn}\ena where
$b(r,\rho;t)$ are transition probabilities of the associated  time
homogeneous birth and death process $\vert X_N^{(\rho)}(t)\vert, \
t\ge 0$ with rates \refm[cor4]. In partucular, if rates of
singular coagulations are positive, then $\rho_r$ is given by
\refm[Q],\refm[Bell] and \refm[ak], while in the case of pure
fragmentation  $\rho_r, \ r=1,\ldots,N$  satisfying \refm[eq2] are as before, if
and only if \refm[fragm] holds.\end{thm}

Note that under $b=0$ in \refm[cor4], the corresponding birth and
death process is known as the   Ehrenfest process (=urn model).\\

\begin{remark} {:\bf Initial distributions} $\rho.$\\
 $CFP's(N)$ with single transitions \refm[psi],\refm[phi]
but  with  initial distributions $\rho$ that do not obey the
equations \refm[eq1], \refm[eq2] in Proposition 2 are not solvable, since in
this case the conditional probability $Q $ depends on $t\ge 0$ and
$\rho,$
 though the processes $\vert X_N^{(\rho)}(t)\vert, \ t\ge 0$ are
time homogeneous.
\end{remark}
\begin{remark} {:\bf  Transition rule for Gibbs fragmentation.}\\
We now explain that the probabilities $P_F$ induced by the
positive fragmentation rates \refm[fragm] define the following
simple rule of a state transition via a fragmentation from
$\eta\in \Omega_{N,r}$ to $\eta_{(i,j)}\in \Omega_{N,r+1}.$ By
\refm[rec], \be P_F(\eta\rightarrow
\eta_{(i,j)})=\frac{(i+j-1)n_{i+j}}{N-r}\
  \begin{cases}
   a_ia_j\
\Big(\frac{1}{2}\sum_{l+m=\ i+j}\ a_la_m\Big)^{-1} , & \text{if}\ i\neq j \\
    \frac{a_i^2}{2} \Big(
\frac{1}{2}\sum_{l+m=2i}\ a_la_m\Big)^{-1}, & \text{if}\ i=j.
  \end{cases}
\la{fragm4}
\end{equation}

Under a given $\eta=(n_1,\ldots, n_N)\in \Omega_{N,r},$ the first
factor in the above expression is the probability that a component
of size $i+j\ge 2$ is selected to fragmentate, while the second
factor specifies the probability that, conditioned on the first
event, the selected component splits into two components of given
sizes $i$ and $j$. As a result, \refm[fragm4] conforms to a
transition procedure postulated in \cite{BerP}, in which the first
and the second  factors are called the linear selection rule and
the Gibbs splitting rule respectively. Theorem  1 and Proposition
2 say that the mean-field transition mechanism \refm[fragm4] is
forced by the requirement that  the process $\vert X_N(t)\vert, \
t\ge 0$ is time homogeneous and the rates of singular coagulations
are positive.
\end{remark}

 \nin {\bf A
historical note} This note concerns exclusively the research on
solvable $CFP's.$ Time evolution of the  stochastic model of pure
coagulation was formulated by Marcus (\cite{mar}) who was also
apparently the first to reveal the relationship between
  Kolmogorov equations and its deterministic analog
 presented by  Smoluchowski equations.
  Solutions to Smoluchowski
 equations for pure coagulation with kernels $K$ induced by
 $\phi(i,j)\equiv const,\ \phi(i,j)=i+j$ and $\phi(i,j)=ij$ were obtained
long ago by researchers in the field of colloid aerosol chemistry
(for references see \cite{mar},\cite{ald1}). Lushnikov (\cite{Lu1}
derived explicit formulae for the expected numbers $\E n_j(t), \
t\ge 0, \  j=1,\ldots ,N$ for the process $X_N^{(\zeta)}(t), \
t\ge 0$ of pure coagulation with $\phi(i,j)=i+j, \ i,j\ge 1,$ with
the help of the generating function for  transition probabilities
$p_{\zeta}(\eta;t), \ t\ge 0,\ \ \zeta,\eta\in \om.$ The
aforementioned stochastic model is known as the Marcus-Lushnikov
process. In (\cite{Lu1}, treating Smoluchowski equations as an
approximation to Kolmogorov ones, Luchnikov proved the important
fact that the solution to Smoluchowski coagulation equations with
a general coagulation kernel,  can be presented as a mixture of
Poisson distributions with time dependent parameters. (Note that
these parameters were found explicitly for  the Marcus -Lushnikov
model only). A further  important contribution was made by
Hendriks, Spouge, Eibl and Schreckenberg (\cite{HendSES}) who
found explicitly the transition probabilities $p_{\zeta}(\eta;t),
\ t\ge 0,\  \zeta,\eta\in \om$ for a more general  Marcus
-Lushnikov model with   $\phi(i,j)$ as in \refm[psi].  This
result, proven via a combinatorial argument, is based on the
representation \refm[prob1] with time independent conditional
probability $Q.$

 \subsection{Discussion of the main result}

\begin{itemize}
\item {\bf Steady state distribution.}
 Firstly, consider solvable $CFP's$ with nonzero rates of  single coagulations and
fragmentations. By \refm[3_8], \refm[FR], \refm[fragm] and Theorem
1, the implied rates $\psi, \phi$ of single coagulations and
fragmentations respectively, are
$$\psi(i,j)=a(i+j)+b,\quad
\phi(i,j)=\phi(1,1)\frac{a_ia_j}{a_{i+j}}\Big(a(i+j)+b\Big), \quad
i,j\ge 1,\quad \phi(1,1)>0,$$  where $a_j, \ j\ge 1$ are given by
\refm[ak]. Thus, the ratio of the above rates is equal to

\be
\frac{\psi(i,j)}{\phi(i,j)}=\frac{a_ia_j}{\phi(1,1)a_{i+j}},\quad
\quad i,j\ge 1.\la{rat12}\end{equation}.

 Setting  in \refm[rat12] \
$\tilde{a}_i=\frac{a_i}{\phi(1,1)},$ gives
$$\frac{\psi(i,j)}{\phi(i,j)}=
\frac{\tilde{a}_i\tilde{a}_j}{\tilde{a}_{i+j}},\quad i,j\ge 1,
$$
 which shows that the criteria of reversibility of mean- field
$CFP's$ (see \cite{DGG})   is fulfilled. Moreover, by Theorem 1,
the above process is the only reversible process, within the class
of solvable mean-field $CFP's.$  By virtue of \refm[dyn], the
invariant measure $\nu_N$ of the process considered is
\be\nu_N(\eta)=b(r)\rho_r(\eta),
 \quad \eta\in
\omr,\ r=1,\ldots,N,  \la{inv}\end{equation} where
$b(r)=\lim_{t\to \infty}b(r,\rho;t)$ is the invariant measure of
the associated  ergodic birth and death process.    The
probability measures $\rho_r$ on $\omr, \ r=1,\ldots,N$ defined by
\refm[Q],\refm[Bell],\refm[ak], belong to the class of
multiplicative measures (= Gibbs distributions) which play also a
role in the theory of random combinatorial structures (see
\cite{V1},\cite{Pi}, \cite{EG}- \cite{frgr2}). The explicit
expression for the measure $\nu_N$ is obtained in a
straightforward way from the known form (see \cite{DGG}) of  the
invariant measure of a reversible $CFP$ with rates \refm[rat12]:
\be \la{cn}\nu_N(\eta)=(c_{N})^{-1}\Big( \frac{\tilde{a}_1^{n_1}
\tilde{a}_2^{n_2} \ldots \tilde{a}_N^{n_N}} {n_1!n_2!  \ldots
n_N!}\Big), \quad \eta\in\om, \end{equation} where $c_{N}$ is the
partition function of the measure $\nu_N.$
 Next, we
embark on an analysis of the asymptotic behaviour of the measure
$\nu_N, $ as $N\to\infty$. For this purpose we need to know the
asymptotics of the weights $\{a_k\}_1^\infty.$ By \refm[ak], \ber
a_k=
  \begin{cases}
   \frac{(\frac{b}{2})^{k-1}}{k!}\prod_{r=2}^k
(\frac{2a}{b}k+r)=
(\frac{b}{2})^{k-1}\big(k!(\frac{2a}{b}k)(\frac{2a}{b}k+1)\big)^{-1}
\Big(\frac{2a}{b}k\Big)_{k+1},   & \ \text{if }\ b\neq 0 \\
\frac{a^{k-1}k^{k-1}}{k!},  & \ \text{if}\ b=0,
  \end{cases} \la{akj}
\ena where $k\ge 1$ and
$(z)_n:=z(z+1)\ldots(z+n-1)=\frac{\Gamma(z+n)}{\Gamma(z)}$ is the
Pochhammer symbol. Applying the Stirling's approximation, gives,
as $k\to \infty,$

\be a_k\sim
  \begin{cases}
C_1C_2^k \ k^{-\frac{3}{2}}     , & \text{if}\ ab>0 \\
    (\frac{b}{2})^{k-1} , & \text{if} \ a=0,\ b>0\\
C_3C_4^k\ k^{-\frac{3}{2}},& \text{if} \ b=0,\ a>0,
  \end{cases}
\end{equation}
where $C_1=C_1(a,b), \ C_2=C_2(a,b), \ C_3=C_3(a), \ C_4=C_4(a)$
are positive constants.  The measure $\nu_N$ in \refm[cn] is
invariant under the transformations of the weights $a_k\rightarrow
C^k a_k,$ with any constant $C>0.$ Thus, the asymptotic behaviour
of the measure $\nu_N$ considered is identical to the one with the
weights \be a_k^\prime\sim
  \begin{cases}
 C_1 \ k^{-\frac{3}{2}}     , & \text{if}\ ab>0 \\
    const , & \text{if} \ a=0,\ b>0\\
C_3\ k^{-\frac{3}{2}},& \text{if} \ b=0,\ a>0,
  \end{cases}
  \la{case}
\end{equation}
as $k\to \infty.$ In accordance with the classification suggested
in \cite{bar} for multiplicative measures $\nu_N$ with regularly
varying weights $a_k\sim k^\alpha, \ k\to \infty,$ the measure
$\nu_N$ considered belongs to the convergent class ($\alpha<-1$)
in the first and the third cases in \refm[case], while in  the
second case in \refm[case] it belongs to the expansive class
($\alpha>-1$). It was shown in \cite{bar} that the convergent
class of $\nu_N$ exhibits a strong gelation, as $N\to \infty$:
with a positive probability all groups cluster in one huge
component of size close to $N.$ In contrast to this, (see
\cite{frgr2}), the expansive measures $\nu_N$ have, with
probability $1,$ as $N\to \infty,$ a threshold value
$N^{\frac{1}{\alpha+2}}$ for the size of the largest group in the
associated  random  partition. In the context of the $CFP$
considered the above crucial difference is easily explained by
noting that the first and the third cases in \refm[case]
correspond to a "strong" coagulation, while the second case
corresponds to a coagulation with a constant rate.

Clearly, pure coagulation and  pure fragmentation processes
$X_N(t), \ t\ge 0$ have the absorbing states $\eta^*=(0,\ldots,1)$
and $\zeta^*=(N,0,\ldots,0)$ respectively.

In the conclusion consider a non- solvable $CFP(N)$ as in Remark
3. In view of the ergodicity of this process, its  invariant
measure will be identical to the one of a solvable $CFP$, starting
from any distribution $\rho$ on $\om$ with the gibbsian
projections $\rho_r$ on $\omr,$ $ r=1,\ldots,N,$ given by
\refm[Q].

\item {\bf CFP's on set partitions}.

These are processes with values in the  space $\Omega_{[N]}$ of
partitions of the set $[N]=\{1,2,\ldots,N\}$ (=set partitions).
From the physical point of view, this means that, in the setting
of this paper, the $N$ particles are labelled, so that clusters
forming a state of the process are subsets of the set $[N].$
$\Omega_{[N]}$- valued $CFP's$ are a generalization  of Kingman's
coalescent that provided a mathematical framework for a variety of
genetic models, in particular the Ewens sampling formula.
Kingman's theory, which is surveyed in \cite{Pi}, is based on the
theory of exchangeable partitions. The development of Kingman's
coalescent by Pitman \cite{Pi} and his colleagues lead to  Gibbs
partitions as distributions of $\Omega_{[N]}$- valued irreversible
processes of pure fragmentation or pure coagulation.  Formally,
the linkage between $\om$- valued and $\Omega_{[N]}$- valued
$CFP's$ is expressed via a simple combinatorial formula and it is
discussed in \cite{Pi}, \cite{BerP}, \cite{EG1}. Among $CFP's$ on
$\Omega_{[N]},$ Gibbs fragmentation processes introduced in
\cite{BerP} by N. Berestycki and Pitman play a central role. These
processes are defined as time homogeneous Markov chains $\Pi(t)\in
\Omega_{[N]}, \ t\ge 0$ of pure fragmentation, such that the conditional
distribution of $\Pi(t)$ given the number of blocks of the random
set partition $\Pi(t)$ is a microcanonical Gibbs distribution not
depending of $t\ge 0$. In terms of CFP's on $\om,$ the above
conditional distributions are just the distributions \refm[Q] on
$\Omega_{N,r}, \ 1\le r\le N.$ Correspondingly, the time reversal of the above
process is called Gibbs coagulation. In \cite{BerP} the authors
posed  a problem of characterization of the weights (in their
notation) $\omega_k:=a_kk!$ for which there exist Gibbs
fragmentation processes, and they proved that, under the
assumption that the fragmentation rates are defined by recursive
and selection rules \refm[fragm4], the unique Gibbs distribution
is given by the weights \refm[rec]. In \cite{BerP}, p.393 it was
conjectured that some other, more complicated splitting rules
might be of interest.
 We will demonstrate (see
Proposition~\ref{abc} below) that the aforementioned
characterization of weights is valid for a broad class of
fragmentation rules, that includes the above one in \cite{BerP}.\\
The problem reduces (see Problem 2 in \cite{BerP}) to
the characterization of weights $a_k, \ k\ge 1, $ (not depending of
$N$) and transition probabilities of fragmentations $P_F$ that
satisfy \refm[eq2], :\ber \la{fr} \rho_{N,r+1}(\zeta)&=&
\sum_{\eta\in \ \rho_{N,r}:\
\eta\sim\zeta}\rho_{N,r}(\eta)P_F(\eta\rightarrow \zeta),\quad \\
\zeta\in \ \Omega_{N,r+1},
   &&\non r=1,\ldots,N-1, \ena when $\rho_{N,r}$ is  a Gibbs
measure \refm[Q] on $\Omega_{N,r}.$ Regarding the probabilities
$P_F,$ we assume that they are of the following general form
implied by the mean- field property: \be P_F(\eta\rightarrow
\eta_{(i,j)})=\frac{n_{i+j}\phi(i,j)}{c(\eta)},\quad
\eta=(n_1,\ldots,n_N)\in\om,\quad \la{asu}
\end{equation}
where $\phi(i,j)$ is a symmetric nonnegative function not
depending of $N$ and \\ $c(\eta)=\sum_{1\le i\le j\le
N}n_{i+j}\phi(i,j)$ is the normalizing constant. Clearly,
\refm[fragm4] is a particular case of \refm[asu].\\
For our subsequent considerations it is important to note that in
\refm[Q] all weights $a_k,\ k\ge 1$ should be  positive, due to
the fact that $\frac{a_N}{B_{N,1}}=1,\ N\ge 1.$

 \begin{proposition}
 \label{abc}

Under the assumption \refm[asu],  Gibbs distributions
$\rho_{N,r},\ r=1,\ldots,N$ satisfy  the equation \refm[fr] if and
only if the weights $a_k$ in \refm[Q]  are given by \refm[rec] and
the rates $\phi(i,j)$ of single fragmentations are the same as in
\refm[fragm].
\end{proposition}

{\bf Proof}\  We assume that Gibbs distributions $\rho_{N,r},\
1\le r\le N$ satisfy \refm[fr]. Treating \refm[fr] when $r=1$ and
$$\zeta=(0,\ldots,\overbrace{1}^{i},\ldots,0,\overbrace{1}^{N-i},0,\ldots,0)\in
 \Omega_{N,2}, \ i=1,\ldots,N-1,$$ gives

\be
\Big(\frac{B_{N,2}}{B_{N,1}}\Big)\Big(\frac{a_N}{a_ia_{N-i}}\Big)\Big(\frac{\phi(i,N-i)}{v_N}\Big)=1,
\end{equation}
if $N\neq 2i$ and
$$\Big(\frac{B_{2i,2}}{B_{2i,1}}\Big)\Big(\frac{a_{2i}}{\frac{1}{2}a_i^2}\Big)
\Big(\frac{\phi(i,i)}{v_{2i}}\Big)=1, $$ if $N=2i,$ where in both
cases $v_k>0$ is defined as in \refm[phi]. Since $B_{N,1}=a_N,$ we
have

\be0<\phi(i,N-i)=\frac{v_N}{B_{N,2}}\begin{cases}
    a_ia_{N-i} , & \text{if } N\neq 2i \\
    \frac{1}{2}a_i^2 , & \text{if } N=2i.\la{sfr}
  \end{cases}
 \end{equation}

Secondly, applying \refm[fr] for $r=2$ with $$\zeta\in
\Omega_{N,3}: \zeta({k_1})=\zeta({k_2})=\zeta({N-k_1-k_2})=1,
$$
where $k_1, k_2, N-k_1-k_2$ are distinct positive integers, gives
\be
1=\sum_{i=1}^3\frac{\rho_{N,2}(\eta_i)}{\rho_{N,3}(\zeta)}P_F(\eta_i\rightarrow
\zeta), \la{eqw}
\end{equation}
where $$\eta_1\in \Omega_{N,2}:\eta_1(k_1)=\eta_1(N-k_1)=1,\quad
\eta_2\in
\Omega_{N,2}:\eta_2(k_2)=\eta_2(N-k_2)=1,$$$$\quad\eta_3\in
\Omega_{N,2}:\eta_3(k_1+k_2)=\eta_3(N-k_1-k_2)=1\quad$$ denote the
three states from which it is possible to arrive, via one step
fragmentation,  at the above state $\zeta\in \Omega_{N,3}.$
Substituting in \refm[eqw] the expression \refm[sfr] for $\phi$
and using \refm[rat], we obtain

\ber \non
1&=&\Big(\frac{B_{N,3}}{B_{N,2}}\Big)\Big(\frac{a_{N-k_1}
v_{N-k_1}}{B_{N-k_1,2}(v_{k_1}+ v_{N-k_1})}\ + \frac{a_{N-k_2}
v_{N-k_2}}{B_{N-k_2,2}(v_{k_2}+ v_{N-k_2})}\ +\\&&
\frac{a_{k_1+k_2} v_{k_1+k_2}}{B_{k_1+k_2,2}(v_{k_1+k_2}+
v_{N-k_1-k_2})}\Big). \la{eqw1} \ena

We set now for a given $N\ge
3,$$$f_N(k):=\frac{a_kv_k}{B_{k,2}(v_k+v_{N-k})},\quad 2\le k\le
N-1.
$$
This allows us to rewrite  \refm[eqw1] as \be
f_N(N-k_1)+f_N(N-k_2)+ f_N(k_1+k_2)=C(N),\quad N\ge 3,\la{rel}
\end{equation} where $C=C(N)$ is a constant w.r.t.  $k_1, k_2:$
$N-k_1\ge 2,\ N-k_2\ge 2,\ k_1+k_2\ge 2.$

The solution of \refm[rel] is given by a linear function
$$f_N(k)=A_Nk+B_N>0, \ 2\le k\le N\ge 3$$ and the constant
$C=2A_NN+3B_N,\ N\ge 3,$ where the reals $A_N,B_N: A_N\ge 0, \
2A_N+B_N>0. $ As a result, the following relation is derived
 \be \frac{a_{k} v_{k}}{B_{k,2}(v_{k}+ v_{N-k})}=A_Nk+B_N,\quad 2\le k\le N\ge
 3.
\la{rlt}
 \end{equation}

We will show that \refm[rlt] forces  the weights $a_k,\ k\ge 2$ to
satisfy  \refm[rec]. Let $$0\le
H_k:=\lim\sup_{N\to\infty}(A_Nk+B_N),$$ for any fixed $k\ge2.$
$H_k=\infty$ is impossible because  $v_k>0, \ k\ge 2,$ by
\refm[sfr]. Hence, $H_k\ge 0$ is finite for all $k\ge 2,$ which
implies  \be 0\le A:=\lim\sup_{N\to\infty} A_N<\infty,\quad
B:=\lim\sup_{N\to\infty} B_N<\infty.\la{456}
\end{equation}

Recalling that $v_1=0,$ we apply \refm[rlt] with $N=k+1, \ k\ge 2$
and $N=2k, \ k\ge 2,$ to get

$$
\frac{a_{k} }{B_{k,2}}=A_{k+1}k +B_{k+1},\ k\ge 2$$

and
$$\frac{a_{k} }{2B_{k,2}}=A_{2k}k +B_{2k},\ k\ge 2$$
respectively. In view of \refm[456], the last two relations are in
agreement if and only if $A=B=0,$ so that from \refm[rlt] we
recognize that
$$ \lim_{N\to \infty}v_{N-k}=\infty, \quad k\ge 1.$$ Consequently,
letting $$z:=\lim\sup_{N\to \infty}\frac{v_{N}}{v_{N-1}}\ge 1,$$
and denoting $$\frac{a_{k}v_k }{B_{k,2}}=e_k>0,\ k\ge 2,$$ one
obtains from \refm[rlt]
$$e_k=\lim_{N\to \infty}v_{N-k}(A_Nk+B_N)=z^{-k}(ak+b), \ k\ge 1, $$
where \be 0\le \tilde{a}:=\lim_{N\to \infty}v_NA_N<\infty,\quad \tilde{b}=:\lim_{N\to
\infty}v_NB_N<\infty.\la{tocef}\end{equation} Substituting the expression for $e_k$ into
\refm[rlt] leads to the following relation
$$\frac{z^{-k}(\tilde{a}k+\tilde{b})}{v_k+v_{N-k}}=A_Nk+B_N,\ 1\le k\le
N-1,$$ which implies

\be\frac{z^{-k}(\tilde{a}k+\tilde{b})+z^{-(N-k)}(\tilde{a}(N-k)+\tilde{b})}{v_k+v_{N-k}}=NA_N+2B_N,\quad
1\le k\le N-1. \la{78}\end{equation}
 Supposing $z>1,$ we obtain, in view of \refm[tocef],
\ber \non z^{-k}(\tilde{a}k+\tilde{b})=\lim_{N\to\infty}v_{N-k}(NA_N+2B_N)=
 &&\\ \begin{cases}
    \infty , &\ \text{if}\ \tilde{a}\neq 0 \\
    z^{-k}\Big( \lim_{N\to\infty}(Nv_{N}A_N)+2\tilde{b}\Big), & \text{otherwise}.
  \end{cases}
 \quad k\ge 1. \ena In both cases this leads to contradiction,
 since in the case $\tilde{a}=0,$ we should have $\tilde{b}>0$, by the definition
of $e_k$.
 Hence, $z=1.$ By \refm[78], this means that for a fixed $N,$ the sum
$v_k+v_{N-k} $ does not depend on $k,$ so that $v_k$ is linear in
$k, $ namely $v_k=\phi(1,1)(k-1),$ since $v_1=0, \ v_2=\phi(1,1).$
As a result, \refm[rlt] becomes
$$\frac{a_k(k-1)}{B_{k,2}}=(N-2)A_{N}k+(N-2)B_N=\frac{\tilde{a}k+\tilde{b}}{\phi(1,1)}:=ak+b, \quad 2\le k\le N\ge 3,$$
since the LHS does not depend on $N$. Recalling now the definition
of $B_{k,2},$ gives \refm[rec]. \qed

\begin{remark} (i) In \cite{BerP}, Berestycki and Pitman
characterized the Gibbs solutions of \refm[fr] in  the particular
case of state transitions \refm[fragm4].  Our solution \refm[ak]
derived under a more restricted mean-field assumption \refm[asu]
has the same form as in  \cite{BerP}, but with less freedom on the
constants $a,b$.

(ii) The weights $w_k=a_kk!$ in the form of a finite product of
linear factors appear also as a solution of a quite different
characterization problem.
 Gnedin and Pitman \cite{gnpit}, extending Kerov's
result \cite{ker}, proved that an infinite sequence $\{\Pi_N\},\
N\ge 1$ of Gibbs random partitions of\ \ $[N]$ is exchangeable  if
and only if in \refm[Q] the weights, say $\tilde{w}_k=\tilde{a}_k
k!,$ are of the form $$\tilde{w}_k
=\prod_{l=1}^{k-1}\Big(\tilde{b}l-\tilde{a}\Big), k\ge 2,\ w_1=1,\
 \tilde{b}\ge 0, \ \tilde{a}\le \tilde{b}.$$ In contrast to
 \refm[ak], the linear factors of $\tilde{w}_k$ do not depend $k$.

\end{remark}

 The  first part of the following corollary gives an answer to  Problem
 3 in \cite{BerP}, in the class of mean- field CFP's, while the second part recovers  Proposition 1 in
 the above paper, in the aforementioned  class of models.

 \begin{cor}
For $N$ large enough there do not exist mean-field Gibbs
fragmentation processes on $\Omega_{[N]}, $ with weights
$w_k=(k-1)!, k\ge 1$ and $w_k\equiv const,\ k\ge 1.$
\end{cor}

{\bf Proof} Recalling that $w_k=a_kk!,\ k\ge 1, $ both assertions
follow from \refm[case] which says that the asymptotics, as $k\to
\infty$ of the two types of  weights in question are not of the
form required in Proposition 3.\qed

\begin{remark} In a recent paper \cite{Gold} a non mean-field Gibbs
fragmentation process with weights $w_k=(k-1)!,\ k\ge 1$ was
constructed. The construction based on the Chinese restaurant
model for simulation of uniform random permutation, results in a
Gibbs fragmentation process with  state transitions not obeying
the mean field form \refm[asu].
\end{remark}

 \item{\bf Spectral
gap}

By virtue of \refm[dyn], the spectral gap of the solvable $CFP's$
considered is equal to the one of the Markov time homogeneous
birth and death process $\vert X_N^{(\rho)}(t)\vert, \ t\ge 0$
with the rates of birth $\lambda_r=\lambda_{r,N}=\phi(1,1)(N-r),\
r=1,\ldots,N-1 $ and rates of death
$\mu_r=\mu_{r,N}=\frac{r-1}{2}\Big(2aN+rb\Big)
 ,\quad
r=1,\ldots,N.$ We shall employ Zeifman's method as described in
\cite{GZ}, to find the spectral gap, say $\beta_N,$ of the above
birth and death process. Recalling that $\lambda_N=\mu_1=0,$
consider the $N-1$ quantities

 \be \alpha_r=\alpha_r(\ \vec{\delta}\ ):=
\lambda_r+\mu_{r+1}-\delta_{r+1}\lambda_{r+1}-\frac{\mu_r}{\delta_r},\quad
r=1,\ldots,N-1, \la{ze}\end{equation} where $\vec
{\delta}=\vec{\delta}_N=( \delta_r=\delta_{r,N}>0,\ r=2,\ldots
N-1)$ is a vector of unknowns $\delta_r$. The method states that\\
(i) For any vector $\vec {\delta},$
$$min \{\alpha_r, \ 1\le r\le N-1\} \le \beta_N\le max\{\alpha_r, \ 1\le r\le N-1\}$$
(ii) In the case of an ergodic birth and death process, there
exists a unique vector $\vec{ \delta},$ such that all $N-1$
quantities $\alpha_r$ are equal, so that their common value is
equal to $\beta_N.$

In our case \refm[ze] conforms to \ber \non
\alpha_r&=&\phi(1,1)(N-r)+\frac{r}{2}\Big(2aN+(r+1)b\Big)
-\phi(1,1)(N-r-1)\delta_{r+1}-\\
&&\frac{(r-1)\Big(2aN+rb\Big)}{2\delta_r},\quad r=1,\ldots,N-1.
\la{ze1}\ena Setting in \refm[ze1] $\delta_r=1, \ r= 2,\dots,N-1,$
we obtain
$$\alpha_r=\phi(1,1)+aN+br, \ r=1,\ldots,N-1,$$
from which the following two -sided bound for the $\beta$ is
derived: $$\phi(1,1)+aN+b\le \beta_N \le \phi(1,1)+aN+b(N-1).$$ In
particular, if $b=0$, the preceding relation gives the exact value
of the spectral gap $\beta_N=\phi(1,1)+ aN.$
\end{itemize}
\vskip .5cm {\bf Acknowledgement}\\ The paper benefitted from helpful remarks and
constructive criticism of three referees. B.G. appreciates an
illuminating discussion with Prof. Aleksander Zeifman.

\end{document}